%% file: ms.tex
\documentclass[reqno,oneside]{amsart}
\usepackage{amsmath,amssymb,amsthm}
\numberwithin{equation}{section}
\PassOptionsToPackage{usenames,svgnames,dvipsnames}{xcolor}
\usepackage{xcolor}
\usepackage[colorlinks=true]{hyperref}
\hypersetup{urlcolor=RubineRed,linkcolor=RoyalBlue,citecolor=ForestGreen}
\usepackage{graphicx}
\newcommand{\PMSp}{P\overline{\mathfrak M}_g^{\mathrm{Sp}(2n)}}
\newcommand{\PMGL}{P\overline{\mathfrak M}_g^{\mathrm{GL}(2n)}}
\newcommand{\Mbar}{\overline{\mathbf M}}
\newcommand{\ThePic}{\mathrm{Pic}(\PMSp)}
\newcommand{\Mgbar}{\overline{\mathcal M}_g}
\newcommand{\Mg}{\mathcal{M}_g}
\newcommand{\CC}{\mathbb{C}}

\newcommand{\eps}{\epsilon}
\newcommand{\vocab}[1]{\emph{#1}}

\newtheorem{theorem}{Theorem}

\begin{document}
\title{Class of discriminant for $\mathrm{Sp}(2n)$ Hitchin spectral covers}
\author{Michael L.\ Baker}

\begin{abstract}
Extending the work \cite{Kor5}, we define the universal Hitchin discriminants in the case of $\mathrm{Sp}(2n)$ Hitchin spectral covers, describe their components, and express their divisor classes in terms of standard generators of the rational Picard group of the moduli spaces of such spectral covers with variable base.
\end{abstract}

\maketitle

\input{introd}
\input{covers}
\input{components}
\input{disc2}
\input{concl}
\input{ack}

\newpage
\input{biblio}

\end{document}

%% file: introd.tex
\section{Introduction}

Hitchin's integrable systems \cite{Hitch1,Hitch2} arise via dimensional reduction of the self-dual Yang--Mills equations. The Hamiltonians of a Hitchin system are encoded in the so-called \emph{spectral cover} $\widehat{\Sigma}$, an $n$-sheeted branched cover of a smooth (or more generally stable) complex projective curve $\Sigma$ defined as a subvariety of $T^* \Sigma$:
\[ \widehat{\Sigma} = \{ (x,v) \in T^* \Sigma \mid P(x,v) = 0 \}, \]
where
\[ P(x,v) = v^n + Q_1(x) v^{n-1} + Q_2(x) v^{n-2} + \cdots + Q_{n-1}(x) v + Q_n(x), \]
with each $Q_j$ a holomorphic $j$-differential on $\Sigma$. Here $v$ can be viewed as a holomorphic 1-form on $\widehat{\Sigma}$. We refer to $P$ (or equivalently, the sequence $(Q_j)$ of its coefficients) as \emph{spectral data}. We denote by $\pi : \widehat{\Sigma} \twoheadrightarrow \Sigma$ the restriction of the canonical projection $T^* \Sigma \twoheadrightarrow \Sigma$ given by $(x,v) \mapsto x$. When all $Q_j$ are arbitrary, this is the notion of a $\mathrm{GL}(n) = \mathrm{GL}(n,\CC)$ spectral cover, studied in \cite{Kor5}. In the framework of \cite{Hitch1,Hitch2} such $P(v)$ arise as characteristic polynomials $P(v) = \det(\Phi - vI)$ of the Higgs field $\Phi$ in the $\mathrm{GL}(n)$ Hitchin system, hence the terminology.

The case of interest for us, however, will be $\mathrm{Sp}(2n) = \mathrm{Sp}(2n,\CC)$ spectral covers: here $P$ is of even degree $2n$ and contains only even powers of $v$. In other words,
\begin{equation}
\label{eqn:spectraldata} 
P(v) = v^{2n} + Q_2(x) v^{2n-2} + \cdots + Q_{2n-2}(x)v^2 + Q_{2n}(x) = \tilde{P}(v^2)
\end{equation}
where
\[ \tilde{P}(q) = q^n + Q_2(x) q^{n-1} + \cdots + Q_{2n-2}(x) q + Q_{2n}(x). \]
(\ref{eqn:spectraldata}) is precisely the form taken by the characteristic polynomial of a Hamiltonian matrix $X$, that is, an element $X \in \mathfrak{sp}(2n)$.

We will be brief, since the ingredients and notation are mostly the same as in \cite{Kor5}. Consider the moduli space $\PMSp$ whose points parametrize the pairs $(\Sigma, [P])$ where $\Sigma$ is a genus $g$ (stable) curve and $P$ is a polynomial of the form (\ref{eqn:spectraldata}) considered up to multiplication by a nonzero constant $\xi \in \mathbb{C}$ given by \[ (\xi \cdot P)(x,v) = \xi^{2n} P(x,\xi^{-1} v). \] There is a natural forgetful map $h : \PMSp \to \Mgbar$ to the Deligne--Mumford compactification of the moduli space of genus $g$ curves, which exhibits $\PMSp$ as a bundle, with weighted projective spaces as fibers. Furthermore, $\PMSp$ is of course a subset of $\PMGL$ (the space denoted $P\overline{\mathcal M}$ in \cite{Kor5} and denoted $P \overline{\mathfrak M}_g^{(2n)}$ in \cite{Bas}), albeit one of rather high codimension. Pulling back the Hodge class $\lambda \in \mathrm{Pic}(\Mgbar)$ along the map $h$, we obtain a class on $\PMSp$ which we also denote by $\lambda$. In this paper $\mathrm{Pic}(\cdot)$ always denotes the \emph{rational} Picard group (that is, coefficients in $\mathbb{Q}$, rather than $\mathbb{Z}$).

When $n=1$, then because $\mathrm{Sp}(2) = \mathrm{SL}(2)$, we see that $\PMSp = P\overline{\mathfrak M}_g^{\mathrm{Sp}(2)}$ is nothing more than a closure of the moduli space of quadratic differentials (considered up to a constant) on genus $g$ smooth projective curves, also known as the projective bundle associated to $T^* \Mg$.

\input{intro/disc}
\input{intro/results}

The structure of the paper is as follows. In \S2 we introduce the various moduli spaces of spectral covers and define the discriminant loci. In \S3 we describe the components of the discriminant loci. In \S4 we prove our main result, Theorem~\ref{fineresult}; the weaker Theorem~\ref{coarseresult} can either be deduced as a consequence (via a routine calculation), or proved independently in an identical fashion. In Section 5 we conclude by listing some interesting open problems.

%% file: intro/disc.tex
\subsection{Discriminant}
\label{discc}

The branch points of the cover $\widehat{\Sigma} \to \Sigma$ are the zeros of the discriminant $W$ of $P$, which is a holomorphic $N$-differential on $\Sigma$ where $N=2n(2n-1)$. It follows that the total number of branch points, counted with multiplicity, is 
\[ 2N(g-1) = 4n(2n-1)(g-1). \]
Standard results on discriminants and their behaviour under polynomial composition imply that, due to the special form (\ref{eqn:spectraldata}) taken by $P(v)$, we have the following factorization:
\begin{equation} 
\label{eqn:discfact}
W = -4^n \cdot Q_{2n} \cdot \Delta^2,
\end{equation}
where $\Delta$ is the holomorphic $2n(n-1)$-differential defined by $\Delta = \mathrm{Discr}(\tilde{P})$. We further define the \emph{reduced discriminant}, a holomorphic $2n^2$-differential, by
\[ W' = -4^n \cdot Q_{2n} \cdot \Delta, \]
so that $W = W' \cdot \Delta$. Generically, $\Delta$ and $Q_{2n}$ have no common zeros, and their zeros are all simple, so $W'$ will have all simple zeros, while the full discriminant $W$ will have double zeros at the former $4n(n-1)(g-1)$ points, and simple zeros at the latter $4n(g-1)$ points. This makes a total of $r=4n^2(g-1)$ zeros, which we denote $x_1, \ldots, x_r$. From the Riemann--Hurwitz formula, we see the genus $\widehat{g}$ of $\widehat{\Sigma}$ in this generic case is
\[ \widehat{g}_{\text{gen}} = (2n)^2(g-1) + 1. \]
Notably, this matches the genus of a generic $\mathrm{GL}(2n,\CC)$-spectral cover. 

The prescription $(\Sigma, [P]) \mapsto (\Sigma, [W])$ defines a map from $\PMSp$ to the moduli space $P\Mbar_g^N$ of pairs $(\Sigma, [W])$ with $[W]$ an $N$-differential on $\Sigma$ considered up to a constant. We of course obtain similar maps using $Q_{2n}$, $\Delta$ and $W'$ in place of $W$; since these differentials generically have all simple zeros, these latter three maps have the additional convenient property that they take a generic point of $\PMSp$ into the principal stratum of the codomain. For this reason the formulas obtained in (the analysis leading to) Theorem~\ref{fineresult} are essentially identical to those of \cite{Kor5}. We note that the loci defined in Section 2 below could also have been naturally defined using these maps.

%% file: intro/results.tex
\subsection{Results}
In \cite{KK} the authors introduced an object called the Bergman tau function (a kind of higher-genus generalization of the Dedekind eta function) on the moduli space of holomorphic Abelian differentials. This moduli space has been extensively studied, due in part to the close relation between Abelian differentials and so-called flat surfaces. In \cite{Kor2}, relations in the Picard group of this moduli space were deduced by studying the asymptotics of the tau function near its boundary. This analytic technique was extended to holomorphic quadratic differentials in \cite{Kor3,Kor3a}, then to certain meromorphic differentials in \cite{KorKall}, then to holomorphic $N$-differentials in \cite{Kor4}, and most recently to $\mathrm{GL}(n)$ Hitchin spectral covers \cite{Kor5,Kor5B}. For a survey of these and many other applications of the Bergman tau function, see \cite{Kor6}. The results of \cite{Kor5} were refined in \cite{Bas}, which in place of the tau function technique uses methods of algebraic geometry (Grothendieck--Riemann--Roch theorem).

In our notation, Theorem 3.2 of \cite{Kor5} reads:

\begin{theorem}
The class of the universal Hitchin discriminant $[PD_W]$ is expressed in terms of the standard generators of $\mathrm{Pic}(P\overline{\mathfrak M}_g^{\mathrm{GL}(n)})$ as follows:
\[ [PD_W] = n(n-1) \Big( (n^2-n+1)(12\lambda - \delta) - 2(g-1)(2n^2-2n+1) \phi \Big). \]
\end{theorem}

Here $\delta = \sum_{j=0}^{[g/2]} \delta_j$ is the pullback of the class of the Deligne--Mumford boundary of $\overline{\mathcal M}_g$, and $\phi$ denotes the first Chern class $c_1(L)$ of the tautological\footnote{The author admits, however, that he still does not know how exactly to associate a tautological line bundle to a \emph{weighted} projective space.} line bundle $L$ on $P\overline{\mathfrak M}_g^{\mathrm{GL}(n)}$ arising from the projectivization.

The most naive approach to $\PMSp$, namely considering the asymptotics of $\tau(\Sigma,W)$ exactly as before, leads to the following:

\begin{theorem}
\label{coarseresult}
Let $n \geq 1$, $g \geq 2$. Putting $N=2n(2n-1)$, the following relation holds in $\ThePic$:
\begin{equation}
\label{eqn:coarseresult}
\lambda = N \kappa_B \phi + \sum_{i=1}^3 c_i [PD_{W,i}] + \frac{1}{12} \delta.
\end{equation}
Here $\kappa_B$ is the homogeneity degree of the Bergman tau function $\tau(\Sigma,W)$ on the stratum $\overline{\mathbf M}_g^N[1^{4n(g-1)},2^{4n(n-1)(g-1)}] = \overline{\mathbf M}_g^N[\mathbf{m}]$ consisting of $N$-differentials with $4n(g-1)$ simple zeros and $4n(n-1)(g-1)$ double zeros. It can be computed as
\begin{equation}
\begin{aligned} 
\label{eqn:coarsehomog}
\kappa_B & = \frac{1}{12N^2} \sum_{k=1}^r \frac{m_k(m_k+2N)}{m_k+N} \\ & = \frac{16n^4-16n^3+12n^2-3n+1}{192n^6 - 288n^5 + 288n^4 - 168 n^3 + 60n^2 - 12n}(g-1) \\
& = \frac{4N^2 + 8N + \sqrt{4N+1} + 5}{12N(N+1)(N+2)}(g-1)
\end{aligned}
\end{equation}
noting that $n = \frac{1}{4}(1 + \sqrt{4N+1})$. Finally,
\begin{align}
\label{eqn:coarsecoeffs}
c_1 & = \frac{1}{12N(N+1)}, \\
c_2 & = \frac{2N+3}{12N(N+1)(N+2)} = \frac{1}{6(N+1)(N+2)} + \frac{1}{4N(N+1)(N+2)}, \\
c_3 & = \frac{1}{3N(N+2)}.
\end{align}
\end{theorem}

(Recall that the relevant tautological line bundles are related by $\mathcal{L} \cong L^N$ and so $\psi = c_1(\mathcal{L}) = N c_1(L) = N \phi$, as in \cite{Kor5}).

It would of course be far more satisfactory to obtain three formulas in which each $[PD_{W,i}]$ occurs alone. This turns out to be possible by examining instead the tau functions $\tau(\Sigma,Q_{2n})$, $\tau(\Sigma,W')$ and $\tau(\Sigma,\Delta)$:

\begin{theorem}
\label{fineresult}
Let $n \geq 1$, $g \geq 2$. In $\ThePic$ the following relations hold:
\begin{equation}
\label{eqn:hodges}
\begin{aligned}
\lambda & = \frac{1}{12 \cdot 2n(2n+1)} [PD_{W,1}] + \frac{(g-1)(4n+1)}{6(2n+1)} \phi + \frac{1}{12} \delta \\
\lambda & = \frac{1}{12 \cdot 8n^2(n-1)} [PD_{W,2}] + \frac{g-1}{3} \phi + \frac{1}{12} \delta \\
\lambda & = \frac{1}{12(2n^2-2n)(2n^2-2n+1)} [PD_{W,3}] + \frac{(g-1)(4n^2-4n+1)}{6(2n^2-2n+1)} \phi + \frac{1}{12} \delta. 
\end{aligned}
\end{equation}
That is, we have the following expressions for the Hitchin discriminants:
\begin{align*}
[PD_{W,1}] & = 2n \Big( (2n+1)(12 \lambda - \delta) - 2(g-1)(4n+1) \phi \Big), \\
[PD_{W,2}] & = 8n^2(n-1) \Big( (12 \lambda - \delta) - 4(g-1) \phi \Big), \\
[PD_{W,3}] & = (2n^2-2n) \Big( (2n^2-2n+1)(12 \lambda -\delta) - 2(g-1)(4n^2-4n+1) \phi \Big).
\end{align*}
\end{theorem}

%% file: covers.tex
\section{Spaces of covers}

We introduce spaces of covers by imitating the $\mathrm{GL}(n)$ case in \cite{Kor5}.

\subsection{Space of covers with fixed base}

Let $\Sigma$ be a smooth curve of genus $g$ and denote by $\mathfrak{M}_{\Sigma}^{\mathrm{Sp}(2n)}$ the moduli space of $\mathrm{Sp}(2n)$ spectral covers of the above form, that is\footnote{In fact we can make this definition for any reductive group $G$ by replacing the $2j$'s with $d_j$, where $d_1,\ldots,d_k$ are the degrees of $G$-invariant polynomials on the Lie algebra $\mathfrak{g}$. As is well-known, for $G$ semisimple we have $\sum (2d_j-1) = \dim G$. The integer $k$ is an invariant of $\mathfrak g$, as are the $d_j$, which are called \emph{fundamental degrees}.},
\[ \mathfrak{M}_{\Sigma}^{\mathrm{Sp}(2n)} = \bigoplus_{j=1}^n H^0(\Sigma, K_\Sigma^{\otimes 2j}) \]
where $K_\Sigma = T^* \Sigma$ is the canonical line bundle on $\Sigma$. Now, recalling that $\dim H^0(\Sigma, K_\Sigma^{\otimes j}) = (2j-1)(g-1)$ for $j \geq 2$ by the Riemann--Roch theorem, we compute
\[ \dim \mathfrak{M}_{\Sigma}^{\mathrm{Sp}(2n)} = \sum_{j=1}^n (4j-1)(g-1) = (2n+1)n(g-1).  \]
Note $(2n+1)n = \dim \mathrm{Sp}(2n)$; in fact, making the analogous definition for general semisimple $G$, it will be the case that
\begin{equation}
\label{eqn:semisimp}
\dim \mathfrak{M}_{\Sigma}^G = (g-1) \dim G.
\end{equation}

We consider three\footnote{The definition of the Hitchin discriminant in \cite{Kor5} does not naively generalize since it is no longer the case that the zeros of $W$ are generically all simple. However it can be generalized by using the reduced discriminant $W'$.} codimension 1 loci in $\mathfrak{M}_{\Sigma}^{\mathrm{Sp}(2n)}$, which we call the \vocab{Hitchin discriminants}. These are, respectively, the loci of spectral data $(Q_j)$ where
\begin{itemize}
\item $Q_{2n}$ has a repeated zero;
\item $Q_{2n}$ and $\Delta$ share one or more zeros;
\item $\Delta$ has a repeated zero.
\end{itemize}
We denote them by $D_{W,1}^\Sigma$, $D_{W,2}^\Sigma$ and $D_{W,3}^\Sigma$ respectively.

\subsection{Space of covers with variable base}

Let $\Mgbar$ be the Deligne--Mumford compactification of the moduli space of genus $g$ curves, and $\nu : \overline{\mathcal C}_g \to \overline{\mathcal M}_g$ be the universal curve. Let $\omega_g = \omega_{\overline{\mathcal C}_g / \overline{\mathcal M}_g}$ be the relative dualizing sheaf, and put
\[ \overline{\mathfrak M}_g^{\mathrm{Sp}(2n)} = \bigoplus_{j=1}^n \Omega_g^{(2j)}, \qquad \overline{\mathfrak M}_g^G = \bigoplus_{j=1}^k \Omega_g^{(d_j)} \]
where $\Omega_g^{(j)} := R^0 \nu_* \omega_g^{\otimes j}$ is the direct image of the $j$th power of $\omega_g$. There is a natural forgetful map $\overline{\mathfrak M}_g^{\mathrm{Sp}(2n)} \to \overline{\mathcal M}_g$ such that the fiber over $\Sigma \in \mathcal{M}_g$ coincides with $\mathfrak{M}_{\Sigma}^{\mathrm{Sp}(2n)}$ (the fibers over nodal curves are a bit more difficult to describe). Thus for semisimple $G$, and in particular $G = \mathrm{Sp}(2n)$, equation (\ref{eqn:semisimp}) gives
\[ \dim \overline{\mathfrak M}_g^{G} = \dim \mathfrak{M}_\Sigma^{G} + \dim \overline{\mathcal M}_g = (\dim G + 3)(g-1). \]
Denote by $D_{W,1}, D_{W,2}, D_{W,3}$ the closures of the sets $\bigcup_\Sigma D_{W,1}^\Sigma, \bigcup_\Sigma D_{W,2}^{\Sigma}, \bigcup_\Sigma D_{W,3}^{\Sigma}$ where in the unions $\Sigma$ runs over $\mathcal{M}_g \subset \overline{\mathcal M}_g$.

As in the $\mathrm{GL}(n)$ case, there is a natural action of $\mathbb{C}^*$ on $\overline{\mathfrak{M}}_g^{\mathrm{Sp}(2n)}$ that looks like $Q_j \mapsto \eps^{j} Q_{j}$ for $\eps \in \mathbb{C}^*$ and $j=1,\ldots,2n$ (note of course that $Q_j=0$ for odd $j$). This respects the loci $D_{W,i}$ defined above ($i=1,2,3$). After projectivization, they become divisors $PD_{W,i}$ in the space $\PMSp$, which we call the \emph{projectivized moduli space of $\mathrm{Sp}(2n)$-spectral data}.

The main goal of this paper is to express the classes $[PD_{W,1}]$, $[PD_{W,2}]$ and $[PD_{W,3}]$ of these divisors in terms of the standard generators of the rational Picard group $\ThePic$.

%% file: components.tex
\section{Components of discriminant loci}

As above, let $r=4n^2(g-1)$. Following \cite{Kor5} we choose a system of generators $(\{ \alpha_i, \beta_i \}_{i=1}^g, \{ \gamma_j \}_{j=1}^r)$ of the fundamental group $\pi_1(\Sigma \setminus \{ x_j \}_{j=1}^r)$ satisfying the standard relation
\[ \gamma_1 \cdots \gamma_r \prod_{i=1}^g \alpha_i \beta_i \alpha_i^{-1} \beta_i^{-1} = \mathrm{id}. \]
The covering $\widehat{\Sigma}$ defines \cite[\S 1.7]{Land} a group homomorphism $G : \pi_1(\Sigma \setminus \{ x_j \}_{j=1}^r) \to S_{2n}$. Let $s_1 = G(\gamma_1)$ and $s_2 = G(\gamma_2)$. The structure of these permutations depends on the type of the zeros $x_1$ and $x_2$. To be precise, $s_k$ will be a simple permutation (that is, a transposition) when $x_k$ is a zero of $Q_{2n}$, while it will be a permutation of type (2,2) when $x_k$ is a zero of $\Delta$. As $x_2 \to x_1$, the covering $\widehat{\Sigma}$ degenerates to a covering $\widehat{\Sigma}_0$ whose structure depends on the type of the product $s_1s_2$. See \cite[\S 5.2.2]{Land} for the terminology ``Maxwell'' and ``caustic''. 

For each case below, we will give $s_1$ and $s_2$ for a representative example. We assume that the sheets have been numbered in such a way that sheets $2k$ and $2k-1$ are conjugate under (that is, interchanged by) the involution, for $k=1,\ldots,n$. It turns out that the presence of the involution in the $\mathrm{Sp}(2n)$ case heavily restricts not only the possible degeneration behaviours, but also the number of \emph{possible inequivalent generic resolutions} in each of those cases (that is, the coefficients appearing in Section~\ref{sec:compsumm}).

\input{comp/PDW1}
\input{comp/PDW2}
\input{comp/PDW3}
\input{comp/summ}

%% file: comp/PDW1.tex
\subsection{$[PD_{W,1}]$: Two zeros of $Q_{2n}$ merge}

Here $s_1$ and $s_2$ are both transpositions. There is only one possibility, due to the privileged status of the point $0 \in T_x^* \Sigma$.

\subsubsection{$[PD_{W,1}^{(b)}]$: Boundary ($n \geq 1$)}

Here $s_1 = s_2 = (12)$, so that $s_1s_2 = ()$. In this case the branching over the newly created double zero of $Q_{2n}$ remains simple. This is because if $x$ is a zero of $Q_{2n}$ which is not a zero of $\Delta$, then the only ramification point in $\pi^{-1}(x)$ is $0 \in T_x^* \Sigma$, which is fixed by the involution. $\widehat{\Sigma}$ develops a nodal singularity (ordinary double point) at this point. We use a box (double vertical lines) to denote a node, while the black colour indicates the point is $0 \in T_x^* \Sigma$. These two conventions will also be followed below.

\begin{center}
\includegraphics[width=100px]{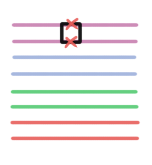}
\end{center}

%% file: comp/PDW2.tex
\subsection{$[PD_{W,2}]$: $Q_{2n}$ and $\Delta$ acquire a common zero}

Here, say $s_1$ is a transposition, while $s_2$ is a permutation of type (2,2). There are two possibilities.

\subsubsection{$[PD_{W,2}^{(ac)}]$: Augmented caustic ($n \geq 2$)}
Here $s_1=(12)$ and $s_2=(13)(24)$, so that $s_1s_2=(1423)$. $\Delta$ is vanishing at $x$ because $\tilde{P}(x,\cdot)$ now has 0 as a double zero (meaning $Q_{2n-2}(x)=Q_{2n}(x)=0$). In this case there still remain $n-2$ nonzero zeros of $\tilde{P}(x,\cdot)$, each giving rise to two distinct points in $\pi^{-1}(x)$, for a total of $2n-4$ points. The only other point in the fibre is $0 \in T_x^* \Sigma$, which is a ramification point of order 4 fixed by the involution, so the fibre is of size $2n-3$, so the total branching order remains 3, so genus $\widehat g$ does not change.

\begin{center}
\includegraphics[width=100px]{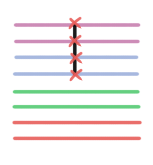}
\end{center}

\subsubsection{$[PD_{W,2}^{(bm)}]$: Boundary--Maxwell ($n \geq 3$)}

Here $s_1=(12)$ and $s_2=(35)(46)$, so $s_1s_2=(12)(35)(46)$. $\tilde{P}(x,\cdot)$ still only has 0 as a simple zero, but two of its \emph{nonzero} zeros have coalesced. This means that $\tilde{P}(x,\cdot)$ only has $n-2$ distinct nonzero zeros, and so as above we see the total branching order remains 3. Three ramification points, each of order 2. Here $0 \in T_x^* \Sigma$ is a ramification point fixed by the involution, and there are two other ramification points which are conjugate under the involution.

\begin{center}
\includegraphics[width=100px]{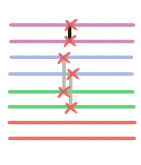}
\end{center}

%% file: comp/PDW3.tex
\subsection{$[PD_{W,3}]$: Two zeros of $\Delta$ merge}

Here, $s_1$ and $s_2$ are both permutations of type (2,2). There are three possibilities, which are just ``pairwise'' versions of the situations from the $\mathrm{GL}(n)$ case.

\subsubsection{$[PD_{W,3}^{(bb)}]$: Doubled boundary ($n \geq 2$)}

Here $s_1 = s_2 = (13)(24)$, so that $s_1s_2 = ()$. $\tilde{P}(x,\cdot)$ still only has one double root. 
Two nodal singularities develop in the fibre above $x$, which are conjugate under the involution. Total branching order is 2, so genus $\widehat{g}$ drops by 1.

\begin{center}
\includegraphics[width=100px]{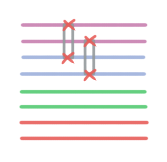}
\end{center}

One might perhaps argue we have missed the case $s_1=(13)(24)$ and $s_2=(14)(23)$, when $s_1s_2=(12)(34)$. This would imply there are two distinct ramification points in $\pi^{-1}(x)$ each fixed by the involution. However, this is impossible: the involution only has one fixed point in $T_x^* \Sigma$, namely $0 \in T_x^* \Sigma$.

\subsubsection{$[PD_{W,3}^{(cc)}]$: Doubled caustic ($n \geq 3$)}

Here $s_1 = (13)(24)$ and $s_2 = (15)(26)$, so that $s_1s_2=(135)(246)$. $\tilde{P}(x,\cdot)$ acquires a root of multiplicity 3. There are two ramification points, each of order 3, conjugate under the involution. Total branching order remains 4, so genus $\widehat{g}$ does not change.

\begin{center}
\includegraphics[width=100px]{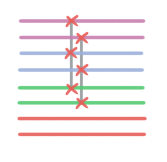}
\end{center}

\subsubsection{$[PD_{W,3}^{(mm)}]$: Doubled Maxwell ($n \geq 4$)}

Here $s_1=(13)(24)$ and $s_2=(57)(68)$, so that $s_1s_2 = (13)(24)(57)(68)$. $\tilde{P}(x,\cdot)$ acquires two double roots. There are four ramification points, each of order 2, which form two conjugate pairs under the involution. Total branching order remains 4, so genus $\widehat{g}$ does not change.

\begin{center}
\includegraphics[width=100px]{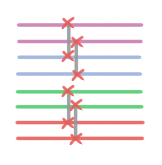}
\end{center}

%% file: comp/summ.tex
\subsection{Summary of components}
\label{sec:compsumm}

We have for $n \geq 4$ that
\begin{align*}
[PD_{W,1}] & = [PD_{W,1}^{(b)}] \\
[PD_{W,2}] & = [PD_{W,2}^{(ac)}] + [PD_{W,2}^{(bm)}] \\
[PD_{W,3}] & = [PD_{W,3}^{(bb)}] + 2 [PD_{W,3}^{(mm)}] + 3 [PD_{W,3}^{(cc)}]
\end{align*}
(for $n < 4$ certain components are absent according to the above).

%% file: disc2.tex
\section{Proof of Theorem~\ref{fineresult}}

We consider three differentials $Q_{2n}$, $W'$, and $\Delta$. The degrees $N_i$ of these differentials are, respectively,
\begin{equation}
\label{eqn:diffdegs}
\begin{aligned}
N_1 & = 2n, \\
N_2 & = N_1 + N_3 = 2n^2, \\
N_3 & = 2n(n-1).
\end{aligned}
\end{equation}
Associated to them are three tau functions $\tau_1 = \tau(\Sigma,Q_{2n})$, $\tau_2 = \tau(\Sigma,W')$, and $\tau_3 = \tau(\Sigma,\Delta)$. The $i$th tau function $\tau_i$ is defined on the principal stratum $\mathbf{M}_g^{N_i}[\mathbf{0}]$ of the moduli space of $N_i$-differentials, which consists of such differentials having all simple zeros. 

The indices $i$ here are consistent with the indices $i$ used to enumerate the divisors $[PD_{W,i}]$, with the exception that unlike $\tau_1$ and $\tau_3$, the tau function $\tau_2$ detects \emph{all three divisors together} rather than $[PD_{W,2}]$ alone. Accordingly, we define $C_i = [PD_{W,i}]$ for $i=1,3$ but $C_2 = \sum_{i=1}^3 [PD_{W,i}]$. Thus $C_i$ is the pullback (along the relevant map; recall the discussion at the end of \S \ref{discc}) to $\PMSp$ of the divisor consisting of $N_i$-differentials with multiple zeros.

Exactly as in \cite{Kor5}, by taking into account the local coordinates near each component of the discriminants $[PD_{W,i}]$, we obtain
\[ C_i = N_i \Big( (N_i+1)(12 \lambda - \delta) - 2(g-1)(2N_i + 1)\phi \Big). \tag{*} \]
For $i=1,2,3$ this reads
\begin{align*}
[PD_{W,1}] & = 2n \Big( (2n+1)(12 \lambda - \delta) - 2(g-1)(4n+1) \phi \Big), \\
\sum_{i=1}^3 [PD_{W,i}] & = 2n^2 \Big( (2n^2+1)(12 \lambda - \delta) - 2(g-1)(4n^2+1)\phi \Big), \\
[PD_{W,3}] & = (2n^2-2n) \Big( (2n^2-2n+1)(12 \lambda -\delta) - 2(g-1)(4n^2-4n+1) \phi \Big).
\end{align*}
Subtracting the first and third from the second, we get
\[ [PD_{W,2}] = 8n^2(n-1) \Big( (12 \lambda - \delta) - 4(g-1) \phi \Big). \]
The right-hand side of this last formula does not fit the pattern of (*). It resembles more the formula for $[P\overline{D}_W^{(m)}]$ or $[P\overline{D}_W^{(c)}]$ in Theorem 2 of \cite{Bas}.

%% file: concl.tex
\section{Conclusion and open problems}

In \cite{Kor5}, a few interesting questions are formulated relating to the $\mathrm{GL}(n)$ case. In this section we propose some further problems.

Firstly, one would of course like to extend the techniques of \cite{Bas} to treat the $\mathrm{Sp}(2n)$ case. This should be fairly straightforward.

Secondly, one would like to carry out this program for an arbitrary complex reductive Lie group $G$. This would presumably involve achieving a deeper understanding of how the structure of the group $G$ controls the various elements involved: the discriminant factorization (\ref{eqn:discfact}), the number of discriminant loci (which turned out to be 1 in the $\mathrm{GL}(n)$ case, but 3 in the $\mathrm{Sp}(2n)$ case), the components of these loci, and so on.

%% file: ack.tex
\section*{Acknowledgements}

The author thanks Dmitry Korotkin and Marco Bertola for discussions.

%% file: biblio.tex
\bibliographystyle{plain}
\addcontentsline{toc}{chapter}{